\input amstex
\documentstyle{amsppt}
\magnification=1200

\vsize19.5cm
\hsize13.5cm
\TagsOnRight
\pageno=1
\baselineskip=15.0pt
\parskip=3pt

\def\p{\partial}

\def\Om{\Omega}

\def\R{\Bbb R}

\def\det{\text{det}}

\def\ol{\overline}
\def\lan{\langle}
\def\ran{\rangle}

\def\phi{\varphi}

\def\osc{\text{osc}}

\nologo
\NoRunningHeads

\topmatter

\title {K\"ahler-Ricci flow on a toric manifold with positive first
Chern class}\endtitle
\author  Xiaohua $\text{Zhu}^{*}$\endauthor
\thanks {\flushpar *  Partially supported by  NSF10425102 in China and a Huo Y-D fund}\endthanks
\address{Department of Mathematics, Peking University,
Beijing, 100871, China}\endaddress
 \email{xhzhu\@math.pku.edu.cn}
\endemail

\abstract  In this note, we prove that on an $n$-dimensional compact
toric manifold  with positive first  Chern class, the K\"ahler-Ricci
flow with any initial $(S^1)^n$-invariant K\"ahler metric  converges
to a K\"ahler-Ricci soliton. In particular, we give another proof
for the existence of K\"ahler-Ricci solitons on a compact toric
manifold with positive first  Chern class by using the
K\"ahler-Ricci flow.
\endabstract

\endtopmatter

\document
\subheading{\bf 0. Introduction}

Let $M$ be  a compact toric manifold with positive first  Chern
class. Let $T\cong (S^1)^n\times \Bbb R^n$ be a maximal torus
which acts on $M$ and $K_0\cong (S^1)^n$ be its maximal compact
subgroup. In this note we discuss a K\"ahler-Ricci flow  with a
$K_0$-invariant initial metric on $M$ and we shall prove

\proclaim {Main Theorem}   On a compact toric manifold $M$ with
positive first  Chern class, the K\"ahler-Ricci flow with any
initial $K_0$-invariant K\"ahler metric  converges to a
K\"ahler-Ricci soliton. In particular, it shows that there exists a
K\"ahler-Ricci soliton on any compact toric manifold  with positive
first  Chern class.
\endproclaim

The existence of K\"ahler-Ricci solitons on a compact toric manifold
with positive first Chern class was  proved in [WZ] by using the
continuity method. The above theorem gives another proof  for the
existence of K\"ahler-Ricci solitons on such a complex manifold by
using the K\"ahler-Ricci flow. We note that a more general
convergence theorem of K\"ahler-Ricci flow on a compact complex
manifold which admits a K\"ahler-Ricci soliton was recently obtained
by  Tian and the author in [TZ3]. In that paper the assumption of
the existence of a K\"ahler-Ricci soliton plays a crucial role. In
the case of K\"ahler-Einstein manifolds with positive first Chern
class the same result was claimed by Perelman ([P2]). In the present
paper we do not need any assumption of the existence of
K\"ahler-Ricci solitons or K\"ahler-Einstein metrics and prove the
the convergence of  K\"ahler-Ricci flow. The Ricci flow was first
introduced by R. Hamilton in 1982 ([Ha]). Recently G. Perelman has
made a major breakthrough in this area for three-dimensional
manifolds ([P1]).

Our proof of the main theorem is to study certain complex
Monge-Amp\`ere flow instead of K\"ahler-Ricci flow. The flow of this
type has been studied before by many people (cf. [Ca], [CT1],
[CT2]). Indeed, our proof used a deep estimate of Perelman ([P2],
also see [ST]). We combined Perelman's estimate with estimates on
solutions of complex Monge-Amp\`ere flow appeared in [TZ3] and used
an argument for $C^0$-estimate on certain real Monge-Amp\`ere
equation studied in [WZ].

The organization of this paper is as follows: In Section 1, We
describe an unpublished estimate of Perelman on the time derivative
of potential functions of evolved K\"ahler metrics along the
K\"ahler-Ricci flow. In Section 2, we reduce  the K\"ahler-Ricci
flow to  a real Monge-Amp\`ere flow in order to get an upper bound
of solutions of potential functions.  Then in Section 3, we use an
argument in [TZ3] to get a $C^0$-estimate of solution. The main
Theorem will be proved in Section 4.

\subheading {\bf 1. An estimate of Perelman}

In this section, we first reduce the K\"ahler-Ricci flow to a
fully nonlinear flow on K\"aher potentials. Then we discuss a
recent and deep estimate of Perelman.

Let $(M,g)$ be an n-dimensional compact K\"ahler manifold with its
K\"ahler form $\omega_g$ representing the first Chern class
$c_1(M)>0$. In local coordinates $z_1,\cdots,z_n$, we have
$$\omega_g=\frac {\sqrt{-1}}{2\pi} \sum_{i,j=1}^n g_{i\overline j}
dz^i\wedge d\overline z^j,~~~~g_{i\bar j}=
g(\frac{\partial}{\partial z_i},\frac{\partial }{\partial z_j}).$$
Moreover, the Ricci form $\text{Ric}(\omega_g)$ is given by
$$\cases &R_{i\overline j}=-\partial_i\overline\partial_j\log
(\text {det}(g_{k\overline l})),\\
&\text{Ric}(\omega_g)=\frac {\sqrt{-1}}{2\pi} \sum_{i,j=1}^n
R_{i\overline j} dz^i\wedge d\overline z^j.\endcases$$ Since the
Ricci form represents $c_1(M)$, there exits a smooth function $h$
on $M$ such that
$$\text {Ric}(\omega_g)-\omega_g=\frac
{\sqrt {-1}}{2\pi}\partial\overline\partial h.\tag 1.1$$

The Ricci flow was first introduced by R. Hamilton in [Ha]. If the
underlying manifold $M$ is  K\"ahler with positive first Chern
class, it is more natural to study the following K\"ahler-Ricci
flow (normalized),

$$\cases &\frac {\partial g(t,\cdot)}{\partial t}
=-\text{Ric}(g(t,\cdot))+g(t,\cdot),\\
&g(0,\cdot)=g_0,\endcases \tag 1.2$$
 where $g_0$ is  a given metric with its K\"ahler class representing $c_1(M)$.
  It can be shown that (1.2) preserves the K\"ahler class,
 so we may write the K\"ahler form of $g(t)$ at a solvable time
$t$ as
$$\omega_{\phi}=\omega_g+\frac
{\sqrt {-1}}{2\pi}\partial\overline\partial \phi$$ for some smooth
function $\phi=\phi(t,\cdot)=\phi_t$. This $\phi$ is usually
called a K\"ahler potential function associated to the K\"ahler
metric $g(t)$. Using the Maximal Principle, one can show that
(1.2) is equivalent to the following complex Monge-Amp\`ere flow
for $\phi(t,\cdot)$,
$$\cases &\frac {\partial\phi}{\partial t}=\log \frac
 {\det (g_{i\overline j}+\phi_{i\overline j})}{\det (g_{i\overline j})}+\phi-h\\
  &\phi (0,\cdot)=0.\endcases\tag 1.3$$
Observe that $\frac{\partial\phi}{\partial t}|_{t=0}=-h$.

Differentiating on both sides of (1.3) on $t$, we have
$$\frac {\partial}{\partial t} \frac {\partial \phi}{\partial t}=
\Delta'\frac {\partial \phi}{\partial t}+\frac {\partial
\phi}{\partial t},$$ where $\Delta'$ denotes the Laplacian
operator associated to the metric $\omega_{\phi}$. Then it follows
from the standard Maximal Principle,
$$|\frac {\partial}{\partial t}\phi(t,\cdot)|\le Ce^t, $$
and consequently,
$$| \phi(t,\cdot)|\le Ce^t.$$
By using these facts and arguments in deriving the higher order
estimates in Yau's solution of the Calabi conjecture [Ya], H.D.
Cao showed that (1.3) is solvable for all $t\in (0,+\infty)$ [Ca].

Using his $W$-functional and arguments in proving non-collapsing
of the Ricci flow [P1], recently, Perelman proved the following
deep estimate [P2] (also see [ST]),

\proclaim {Lemma 1.1} Let $\phi_t$ be a solution of Monge-Amp\`ere
flow (1.3).  Choose $c_t$ by the condition $h_t=-\frac {\partial
\phi}{\partial t}+c_t$ such that
$$\int_M e^{h_t}\omega_{\phi}^n=\int_M\omega_g^n.$$
Then there is a uniform constant $A$ independent of $t$ such that
$$|h_t|\le A.\tag 1.4$$
\endproclaim

Lemma 1.1 is crucial in proving our main theorem. Recall that
$h_t$ is defined by (1.1) with $\omega_g$ replaced by
$\omega_\phi$ and can be different to a constant. In Section 3
below we will further prove that $c_t$ is  uniformly bounded and
so $\frac{\partial \phi}{\partial t}$ is.

\subheading {\bf 2.  Upper bound of solution }

In this section,  we discuss the upper bound of solution of
equation (1.3) by reducing it to a real Monge-Amp\`ere flow. We
now assume that $M$ is a compact toric manifold with positive
first Chern class ( toric Fano manifold) and $g$ is  a
$K_0$-invariant K\"ahler metric on $M$. Then under an affine
logarithm coordinates $(w_1,...,w_n)$, its K\"ahler form
$\omega_g$  is determined by a convex function $u$ on $\R^n$,
namely
$$\omega_{g}=\frac {\sqrt {-1}}{2\pi}
         \p \ol \p  u \ \ \ \text {on }~  T.$$
Hence
$$\omega_{g}^n=(\frac {1}{\pi})^n  \text {det}(u_{ij})
      dx_1\wedge ...\wedge dx_n\wedge d\Theta, $$
where $w_i=x_i+\sqrt{-1}\theta_i$ and $d\Theta=d\theta_1\wedge
...\wedge d\theta_n$ is the standard volume form of $K_0$.

Let $\Om\subset\Bbb R^n$ be a  bounded convex polyhedron
associated to the toric Fano manifold $M$ and denote $p^{(1)},
\cdots, p^{(m)}$ to be the vertices of $\Om$. We define a convex
function  on $\Bbb R^n$ by
$$v^0(x)=\log\big(\sum_{i=1}^m e^{\lan p_i, x\ran}\big). \tag 2.1$$
Then the induced metric $\omega_{g^0}=\frac {\sqrt {-1}}{2\pi}
         \p \ol \p  v^0  $ can be extended as one with $c_1(M)=[\omega_{g_0}]$ on $M$ [BS].
The expression (2.1)  implies that the gradient (moment) mapping
$Dv^0$ is a diffeomorphism from $\Bbb R^n$ to $\Omega$ and
$$|\log \text {det}(v^0_{ij}) + v^0| < \infty .  \tag 2.2$$
Denote
$$\ol v (x)= \max\{x\cdot p^{(k)}\ {|}\ k=1, \cdots, m\}. \tag 2.3$$
The graph of $\ol v$ is a convex cone with vertex at the origin.
It can be verified that
$$ |\ol v-v^0|\le C,\tag 2.4$$
namely the graph of $\ol v$ is an asymptotical cone of the graph
$v^0$.

Let $h_0$ be a smooth function determined by the relation (1.1)
associated to the metric $\omega_{g^0}$. Then it is clear,
$$\p \ol \p  [e^{v^0+h_0}\text{det}(v^0_{ij})]=0
    \ \ \ \text {in}~\Bbb R^n. $$
Hence by (2.2) we have, after normalization,
$$ \text{det}(v^0_{ij})=e^{-h_0-v^0}. \tag 2.5$$
In general, the equation still holds for a convex function $u$ on
$\Bbb R^n$ induced by a $K_0$-invariant K\"ahler metric on $M$
since the difference between $u$ and $v_0$ can be extended as a
smooth function on $M$. Note that $\text{Im}(Du)=\Omega,$ where
$\text{Im}(Du)$ denotes the image of gradient map of $u$ in $\Bbb
R^n$.

Now we consider  K\"ahler-Ricci flow (1.2) with a $K_0$-invariant,
initial K\"ahler metric $g$ which is  induced by a convex function
$u_0$ on $\Bbb R^n$. For simplicity, we may assume that $u_0$
satisfies $\inf u_0= u_0(o)=0$. Since $K_0$-invariant preserves
under the flow, we can reduce equation (1.3) to a real
Monge-Amp\`ere flow as follow,
$$\cases &\frac {\partial u}{\partial t}=\log\text {det} (u_{ij})+u,
              \ \ \ \text{in}~ \Bbb R^n,\\
              &u(0,\cdot)=u_0, \endcases \tag 2.6$$
where $u=u_0+\phi.$

\proclaim{Lemma 2.1} Let $u=u_t=u(t,\cdot)$ be a solution of
equation (2.6) and $\overline u=\overline u_t=u_t-c_t$, where
$c_t$ are functions appeared in Lemma 1.1. Let $m_t=\inf_{\R^n}
\overline u_t(x).$ Then
 $$ |m_t| \le C$$
  for some $C>0$ independent of $t\in (0, \infty)$.
\endproclaim

\demo{ Proof} By Lemma 1.1, we have
 $$|\frac {\partial u}{\partial t}-c_t|\le A.$$
 It follows by equation (2.6),
 $$c_1\le\int_{\Bbb R^n} e^{-\overline u} dx\le c_2,$$
   for some uniform constants $c_1$ and $c_2$. Note that
   $$|D\overline u_t|\le \text{diam}(\Omega),$$
   where $ \text{diam}(\Omega)$ denotes the diameter of $\Omega$.
Then it is easy to see that $m_t$ is uniformly bounded from below.

 To get an upper bound of $m_t$, we
use an argument in [WZ]. For any nonnegative integer $k$, we
denote a set,
$$A_k=\{x\in\R^n:\ m_t+k\le \overline u(x)\le m_t+k+1\}. $$
Then for any $k\ge 0$,  set $\bigcup_{i=0}^k A_i=\{w<m_t+k+1\}$ is
convex. Note that  the origin is contained in $\Om$.  Hence the
minimum $m_t$ is attained at some point in $A_0$ and $A_k$ is a
bounded set for any $k\ge 0$. By a well-known theorem [Mi], there
is a unique ellipsoid $E$, called the minimal ellipsoid of $A_0$,
which attains minimum volume among all ellipsoids contain $A_0$ ,
such that
$$\frac {1}{n}E\subset A_0\subset E.$$
Let $B$ be a linear transformation with $\det (B)=1$, which leaves
the center of $E$ invariant, such that $B(E)$ is a ball $B_R$ with
radius $R$. Then we have $B_{R/n}\subset B(A_0)\subset B_R$ for
two balls with concentrated center.

 By equation (2.6) and Lemma 1.1,  we have
$$\aligned \det  (\overline u_{ij})=& \exp\{\frac{\partial u}{\partial t}-c_t-\overline u \}\\
&\ge c e^{-\overline u}, ~~\text {in}~ \Bbb R^n. \endaligned$$
 It follows
$$\det  (\overline u_{ij})\ge \frac {c}{e} e^{-m_t}, \ \ \ \text{in}\ \ A_0. $$
 We claim
$$R\le \sqrt 2\, n (\frac {c}{e})^{-1/2n}\, e^{m_t/2n}. \tag 2.7$$

 Let
$$v(\overline y)=\frac {1}{2}(\frac {c}{e})^{1/n}e^{-m_t/n}
           \big[|\overline y-\overline y_t|^2-(\frac {R}{n})^2\big] +m_t+1,$$
where $\overline y_t$ is the center of the minimum ellipsoid of
$A_0$. Then
$$\det  (v_{ij})= \frac {c}{e} e^{-m_t},\ \ \ \text{in}~ B(A_0),$$
and $v\ge \overline u$ on $\p B(A_0)$, where $\overline u=\overline
u(B^{-1})$.  Hence by the comparison principle for the
Monge-Amp\'ere operator we have $v\ge \overline u$ in $B(A_0)$. In
particular, we have
$$\align
m_t & \le \overline u(\overline y_t)\le v(\overline y_t)\\
    & = -\frac {1}{2}(\frac{c}{e})^{1/n}e^{-m_t/n}
           (\frac {R}{n})^2+ m_t+1.\\
\endalign $$
Hence (2.7) follows.

By the convexity of $\overline u$, we have
$$B(A_k)\subset B_{2(k+1)R}(0).$$
Thus by (2.7),  we obtain
$$\align
\int_{\R^n}e^{-\overline u} d\overline y& =\sum_k
\int_{B(A_k)}e^{-\overline u}
d\overline y\\
&=\sum_k \int_{B(A_k)}e^{-\overline u}|B(A_k)|\\
        & \le \omega_n \sum e^{-m_t-k} |2(k+1)R|^n\\
    & = \omega_n\frac {(2R)^n}{e^{m_t}}\sum \frac {(k+1)^n}{e^k}\\
    & \le C_1 e^{-m_t/2}, \endalign $$
where $\omega_n$ is the area of the sphere $S^{n-1}$. We notice
that the above integration is invariant under any linear
transformation $B$ with $|B|=1$. Returning to the original
coordinates $x$, by equation (2.6), we have
$$\aligned
e^{-m_t/2} & \ge \frac {1}{C_1} \int_{\R^n} e^{-\overline u}dx\\
 & = \frac {1}{C_1}\int_{\Bbb R^n} \exp\{-\frac {\partial u}{\partial t}+c_t\}\det(\overline u_{ij})dx\\
   &\ge \frac{e^{-A}}{C_1} \int_{\Bbb R^n}\det(\overline u_{ij})dx=\frac{e^{-A}}{C_1}|\Omega|.
\endaligned $$
Hence $m_t\le C$.  \qed \enddemo

Let $x_t$ be the minimal point of $\overline u_t$ and $\overline {
\overline u}=\overline {\overline u}_t(\cdot) =\overline u_t(\cdot
+x_t)-m_t$. Set $\overline\phi=\overline{\overline u}-u_0$. Then

\proclaim {Proposition 2.1}
$$|\sup_M\overline\phi |\le C.$$
\endproclaim

\demo {Proof} The upper bound of $\sup_M\overline\phi$ follows
from the fact,
$$\aligned \overline \phi&=\overline{\overline u}- u_0= \overline{\overline u}- \overline v
+\overline v-u_0\\
&\le \overline v-u_0=|\overline v-u_0|\\
&\le |\overline v- v_0|+|v_0-u_0|\le C.\endaligned$$

 From (2.6), we obtain an equation for $\overline {\overline
u}(.)$,
$$\text{det}(\overline{\overline u}_{ij})= \exp\{\frac{\partial u}{\partial t}
-c_t-\overline{\overline u}-m_t\},
~~\text{in} ~\Bbb R^n.\tag 2.8$$
 It follows that $\overline\phi$ satisfies an equation,
$$\text{det}(g_{ij}+\overline\phi_{ij})=
 \exp\{h+\frac{\partial u}{\partial t}-c_t-\overline\phi-m_t\},
 ~~\text {in}~M,\tag 2.9$$
where $h$ is  the smooth function determined by the relation (2.5)
as $v_0$ is replaced by $u_0$.  By (2.9), it is easy to see that
$$c_1\le\int_Me^{-\overline\phi-m_t}\omega_g^n\le c_2, \tag 2.10$$
 for some positive
constants $c_1$ and $c_2$ since $\frac{\partial u}{\partial
t}-c_t$ are uniformly
 bounded. Thus
$$\sup_M\overline \phi+m_t\ge -C_0~\text{and}~\inf_M\overline \phi+m_t\le
C_0\tag 2.11$$
  for some uniform constant $C_0$. Therefore  by Lemma 2.1, we get
$$\sup_M\overline\phi \ge C.$$
The proof is completed.\qed\enddemo

\subheading {3. Generalized K-energy and $C^0$-estimate}

In this section, we use the monotonicity of generalized K-energy
introduced in [TZ2] and a $C^0$-estimate developed in [TZ3] to get
a bound of $C^0$-norm of $\overline\phi$.

Let $\eta(M)$ be a Lie algebra which consists of all holomorphic
vector fields on $M$ and $\eta_0(M)$ be a Lie subalgebra of
$\eta(M)$ induced by torus $T$.  According to [WZ], we see that
there is  a holomorphic vector field  $X\in \eta_0(M)$ on $M$ such
that the holomorphic invariant $F_X(\cdot)$ introduced in [TZ2]
vanishes , i.e.,
 $$F_X(\cdot)\equiv 0,~\text{on }~\eta(M).$$

    Set a subspace of potential functions space by
$$\Cal M_X(\omega_g) =\{\phi\in C^{\infty}(M)|
~\omega_{\phi}=\omega_g+ \frac {\sqrt {-1}}{2\pi}
\partial\overline\partial\phi >0,~~\text{Im}(X)(\phi)=0\}.$$
 The generalized K-energy associated to $X$ is  a functional
 defined in $\Cal M_X(\omega_g)$ by
$$\aligned & \tilde\mu(\phi)\\
 &=-\frac {n\sqrt {-1}}{2\pi V}\int_0^1\int_M
\dot\psi[\text
{Ric}(\omega_{\psi})-\omega_{\psi}\\
&-\frac{\sqrt{-1}}{2\pi}(\partial\overline\partial
\theta_X(\omega_{\psi}) -\partial
(h_{\omega_{\psi}}-\theta_X(\omega_{\psi}))\wedge\overline\partial\theta_X(\omega_{\psi}))
]\wedge e^{\theta_X(\omega_{\psi})}\omega_{\psi}^{n-1} \wedge
dt,\endaligned$$
 where $\psi=\psi_t$ $(0\le t\le 1)$ is a path connecting $0$ to $\phi$ in $\Cal
M_X(\omega_g)$ and $\dot\psi=\frac{d\psi}{dt}$,  and
$\theta_X(\omega_{\psi})=\theta_X +X(\psi)$ with $\theta_X$
defined by
$$\frac{\sqrt{-1}}{2\pi}\overline\partial\theta_X=i_X(\omega_g).$$
We notice that $\tilde\mu(\phi)$ is just Mabuchi's $K$-energy [Ma]
when $X=0$ in the relation. One can show that for any
$\sigma\in\text{Aut}_r(M)$, $\phi_{\sigma}\in\Cal M_X(\omega_g)$
and
$$\tilde\mu(\phi)=\tilde\mu(\phi_{\sigma})\tag 3.1$$
where $\text{Aut}_r(M)$ denotes a reductive subgroup of the
holomorphic automorphisms group on $M$ containing $T$ and
$\phi_{\sigma}$ is defined by
 $$\sigma^*\omega_{\phi}=\omega_g+ \frac
 {\sqrt {-1}}{2\pi}\partial\overline\partial\phi_{\sigma}.$$

 Let $\sigma_t=\exp\{tX\}$
be the one-parameter group generalized by $X$ and
 $\phi'=\phi'_t$ be  a family of potential functions given by relations
 $$\aligned \sigma_t^*\omega_{\phi_t}&
 =\sigma_t^*\omega_g+ \frac
 {\sqrt {-1}}{2\pi}\partial\overline\partial(\phi_t\cdot\sigma_t)\\
 & =\omega_g+ \frac
 {\sqrt {-1}}{2\pi}\partial\overline\partial\phi'_t.\endaligned$$
  Then  according to equation (1.3), $\phi'$ modula   constants
  satisfy
  the following complex Monge-Amp\`ere flow,
$$\frac {\partial\phi'}{\partial t}=\log \frac
 {\det (g_{i\overline j}+\phi_{i\overline j}')}{\det (g_{i\overline j})}
 +X(\phi')+\phi'-h+\theta_X. \tag 3.2$$
By using this equation,  one can show
$$\aligned &\frac {d\tilde\mu(\phi')}{dt}\\
&=-\frac {1}{V}\int^{\infty}_0\int_M\|\overline\partial\frac
{\partial\phi'}{\partial t}\|^2
 e^{\theta_X+X(\phi')}(\omega_{\phi'})^n.\endaligned\tag 3.3$$
In particular,
$$\tilde\mu(\phi')\le 0.$$
Therefore by (3.1), we get
$$\tilde\mu(\overline\phi)=\tilde\mu(\overline\phi+m_t)
=\tilde\mu(\phi)=\tilde\mu(\phi')\le 0.\tag 3.4$$

 \proclaim {Proposition  3.1}
$$\|\overline\phi\|_{C^0(M)}\le C.$$
\endproclaim

\demo {Proof} Recall another  functional introduced in [TZ2],
$$\aligned &\tilde F(\phi)\\
&= \tilde J(\phi)-\frac{1}{V}\int_M\phi e^{\theta_X}\omega_g^n
-\log(\frac{1}{V}\int_M e^{h-\phi}\omega_g^n),~\text{in}~\Cal
M_X(\omega_g),\endaligned$$
  where
  $$\tilde J(\phi)=\frac
{1}{V}\int_0^1\int_M\dot\psi(e^{\theta_X}\omega_g^n
-e^{\theta_X+X(\psi)}\omega_{\psi}^n)\wedge dt>0$$
 is a functional in $\Cal M_X(\omega_g)$ which is independent of the choice of path
 $\psi_t$ $(0\le t\le 1)$ connecting $0$ and $\phi$ [Zh]. It was proved in [TZ2],
 $$\tilde \mu(\phi)\ge \tilde F(\phi)-C,~\forall ~\phi\in\Cal M_X(\omega_g), \tag 3.5$$
for some uniform constant $C$. Thus by (3.4), we have
$$\tilde F(\overline\phi+m_t)\le C.\tag 3.6$$

  By (2.10) and Proposition 2.1, we get from (3.6),
  $$ \tilde J(\overline\phi)=\tilde J(\overline\phi+m_t)\le \frac{1}{V}\int_M(\overline \phi+m_t)
e^{\theta_X}\omega_g^n+C\le C_1.\tag 3.7$$
 Since
  $$\tilde J(\overline\phi)\ge c I(\overline\phi)$$
for some uniform constant $c$ [CTZ], where
$$I(\overline \phi)=\frac {1}{V}\int_M \overline\phi(\omega_g^n
-\omega_{\overline\phi}^n),$$
 we obtain
 $$I(\overline \phi)\le C_2.\tag 3.8$$
  On the other hand, applying the argument in the proof of Proposition 3.1 in
[TZ3] to equation (2.9), one can show that
$$\text{ocs}_M\overline\phi\le C_0(1+I(\overline \phi))^{n+1}$$
for some uniform constant $C_0$. Then  by (3.8), it follows
$$\text{ocs}_M\overline\phi\le C_3.$$
Thus by (2.11), it is easy to see
 $$\|\overline \phi+m_t\|_{C^0(M)}\le \text{ocs}_M(\overline\phi+m_t)+2\tilde C
 =\text{\osc}_M\overline\phi+2\tilde C\le C_4.$$
Therefore, By Lemma 2.1, we prove  Proposition 3.1. \qed\enddemo

\proclaim {Corollary 3.1} Let $\phi=\phi_t$ be a  solution of
equation (1.3). Then
$$\tilde \mu(\phi)\ge C.\tag 3.9$$
\endproclaim

\demo {Proof} By (2.10) and (3.5), we have
 $$\aligned \tilde \mu(\phi)&=\tilde\mu(\overline\phi+m_t)\\
 &\ge\tilde F(\overline\phi+m_t)-C_1\\
 &\ge -\int_M(\overline \phi+m_t) e^{\theta_X}\omega_g^n-C_2.\endaligned$$
   Then (3.9) follows from Lemma 2.1 and Proposition 3.1.\qed\enddemo

By (3.1), (3.3) and Corollary 3.1,  we have
$$\frac {1}{V}\int^{\infty}_0\int_M\|\overline\partial\frac
{\partial\phi'}{\partial t}\|^2
 e^{\theta_X+X(\phi')}(\omega_{\phi'})^n\wedge dt\le C.$$
It follows
$$\frac {1}{V}\int^{\infty}_0\int_M\|\overline\partial\frac
{\partial\phi'}{\partial t}\|^2
 e^{\theta_X+X(\phi')-t}(\omega_{\phi'})^n\wedge dt\le C.\tag 3.10$$

The following proposition was proved in [TZ3].

  \proclaim {Proposition 3.2} Let $h$ be normalized by adding a suitable
constant so that
$$\frac {1}{V}\int_M (h-\theta_X)\omega_g^n
 =-\frac {1}{V}\int^{\infty}_0\int_M\|\overline\partial\frac
{\partial\phi'}{\partial t}\|^2
 e^{\theta_X+X(\phi')-t}(\omega_{\phi'})^n\wedge dt.\tag 3.11$$
Then
 $$|c_t|\le C~\text{and}~|\frac {\partial\phi}{\partial t}|\le C,$$
where $c_t$ are constants appeared in Lemma 1.1. \endproclaim

\subheading{4. Convergence of the flow}

In this section, we discuss the higher order estimates for a
modified solution of equation (1.3) and finish the proof of Main
Theorem.  Let $x_t\in \Bbb R^n$ be  a family of points determined
in section 2. We observe that

\proclaim {Lemma 4.1} Let $i=0,1...,$ be any nonnegative integer.
Then the distances between $x_i$ and $x_{i+1}$ are uniformly
bounded, i.e.,
$$|x_{i}-x_{i+1}|\le C$$
 for some uniform constant $C$.
 \endproclaim

\demo {Proof} Let $\overline v$ be the convex cone function in
$\Bbb R^n$ defined by (2.3).  Then by Proposition 3.1, we have
$$|\overline {\overline u}_i-\overline v|
\le |\overline {\overline u}_i- u_0|+| u_0-\overline v|=
|\overline\phi_i|+| u_0-\overline v|\le C,$$
 and
 $$|\overline {\overline u}_{i+1}-\overline v|
\le |\overline {\overline u}_{i+1}- u_0|+| u_0-\overline v|=
|\overline\phi_{i+1}|+| u_0-\overline v|\le C.$$
 On the other hand, by Proposition 3.2, we have
$$\aligned &|\overline{\overline u}_i(x-x_i)-\overline{\overline
u}_{i+1}(x-x_{i+1})|\\
&=|\overline u_i(x-x_i)-\overline u_{i+1}(x-x_{i+1})
-m_i+m_{i+1}|\\
&
=|u_i-u_{i+1}-m_i+m_{i+1}|(x)=|\phi_i-\phi_{i+1}-m_i+m_{i+1}|(x)\le
C.\endaligned$$
 Thus combining the above  relations, we get
$$|\overline v(x-x_i)-\overline
v(x-x_{i+1})|\le C.$$
 In particular, by choosing $x=x_i$, we see
$$|\overline v(x_{i}-x_{i+1})|
=|\overline v(o)-\overline v(x_i-x_{i+1})| \le C.\tag 4.1$$
 (4.1) implies
 $$|x_i-x_{i+1}|\le C, ~~\forall ~i=0,1,....$$
 \qed\enddemo

With the help of  Lemma 4.1,   one can choose a  family of
modified points $x'_t$ in $\Bbb R^n$ such that
$$|x_t-x'_t|\le C~\text{and}~|\frac{dx'_t}{dt}|\le C.\tag 4.2$$
Let $\tilde u(\cdot)=\tilde u_t(\cdot)=u_t(x'_t+\cdot)$ and
$\tilde \phi=\tilde\phi_t=\tilde u_t-u_0.$ Then
$\tilde\phi_t(\cdot)=\phi_t(x_t'+\cdot)$ is a potential function
on $M$.  Since
$$|\tilde u_t-\overline{\overline u}_t|\le \text {diam}(\Omega)|x_t-x'_t|+|c_t|+|m_t|,$$
by Proposition 3.1, we have
$$\|\tilde\phi_t\|_{C^0(M)}\le\|\tilde u_t-\overline{\overline u}_t \|_{C^0(M)}+
\|\overline\phi_t\|_{C^0(M)} \le C.$$

From (2.6), we see that  $\tilde u$ satisfies a parabolic  equation
$$\frac{\partial\tilde u}{\partial t}= \log\text{det}(\tilde u_{ij})+\tilde X(\tilde u)+\tilde u,
~~\text{in} ~\Bbb R^n,\tag 4.3$$
 where $\tilde X=\tilde X_t=\frac {dx'_t}{dt}$. Note that  $\tilde
 X$ are corresponding to
 a family of holomorphic vector fields on $M$. It follows that $\tilde\phi$ satisfies an equation
$$\aligned &\frac{\partial\tilde\phi}{\partial t}\\
&= \log\text{det}(g_{ij}+\tilde\phi_{ij})-\log
 \text{det}(g_{ij})+\tilde X(\tilde\phi)+\tilde\phi-h+\theta_{\tilde X_t},
 ~~\text {in}~M,\endaligned\tag 4.4$$
where $h$ is normalized as in Proposition 3.2. Note that
$\theta_{\tilde X_t}=\tilde X_t(u_0)$. So by Proposition 3.2, we
have
$$|\frac{\partial\tilde\phi}{\partial t}|=|\frac{\partial\phi}{\partial t}+\theta_{\tilde X_t}|\le C$$

\demo {Proof of Main Theorem}  We shall show that the
corresponding K\"ahler metrics $\omega_{\tilde\phi}$ associated to
solution of equation (4.4) converge to a K\"ahler-Ricci soliton
$\omega_{KS}$ with respect to $X$. The proof is similar to one of
Main theorem in [TZ3]. We give a sketch. First by modifying Yau's
$C^2$-estimate in [Ya] for certain complex Monge-Amp\`ere
equation, one obtains for the solution $\tilde\phi$ of equation
(4.4),
$$\|\tilde\phi\|_{C^2(M)}\le
C~\text{and}~(g_{ij}+\tilde\phi_{ij})>c_0.$$
 Then following Calabi's $C^3$-estimate for complex Monge-Amp\`ere equation [Ya],
  we further get
$$\|\tilde\phi\|_{C^3(M)}\le C.$$
Thus by using regularity theory for the parabolic equation, one
sees easily that all $C^k$-norms of the solution $\tilde\phi$ are
uniformly bounded. Therefore,  we conclude  that for any sequence
of functions $\tilde\phi_t$, one can take a subsequence of the
sequence which converge  $C^k$-smoothly to a smooth function
$\overline\phi_{\infty}$ on $M$.

Let $\sigma_t=\exp\{tX\}$ and
$\sigma'=\sigma_t'=\rho_t\cdot\sigma_t^{-1}$, where $\rho_t$ are
holomorphic transformations corresponding to changes from $\phi_t$
to $\tilde\phi_t$. Since
$$\int_M\|\overline\partial ((\sigma')^*\frac {\partial\phi'}{\partial t})\|^2
 e^{\theta_X+X(\tilde\phi)}\omega_{\tilde\phi}^n
=\int_M\|\overline\partial \frac {\partial\phi'}{\partial t}\|^2
 e^{\theta_X+X(\phi')}(\omega_{\phi'})^n,$$
then by (3.3), we have
$$ \int_M\|\overline\partial ((\sigma')^*\frac {\partial\phi'}{\partial t})\|^2
 e^{\theta_X+X(\tilde\phi)}\omega_{\tilde\phi}^n
=-\frac {d\tilde\mu_{\omega_g}(\phi')}{dt}.\tag 4.5$$
  By the  lower bound of $\tilde\mu_{\omega_g}(\phi')$ (cf. Corollary 3.1), one sees that
there is a sequence of $t_i$, $i=1,2,...,$ such that
$$\int_M\|\overline\partial ((\sigma')^*\frac {\partial\phi'}{\partial t}|_{t_i})\|^2
\omega_{\tilde\phi_{t_i}}^n\to 0,~\text {as}~i\to\infty.\tag 4.6$$
 On the other hand, by (4.4), we have
$$\frac{\sqrt{-1}}{2\pi}\partial\overline\partial
[(\sigma')^*\frac {\partial\phi'}{\partial t}] =-\text
{Ric}(\omega_{\tilde\phi})+\omega_{\tilde\phi}+L_X\omega_{\tilde\phi}.\tag
4.7$$
 Then $(\sigma')^*\frac {\partial}{\partial t}\phi'$ is
$C^k$ uniformly bounded in space, so  there exists  a convergent
subsequence of $(\sigma')^*\frac {\partial}{\partial
t}\phi'(t_i,\cdot)$. Hence by (4.6), we conclude that
$(\sigma')^*\frac {\partial}{\partial t}\phi'(t_i,\cdot)$ (still
use same indices $t_i$) converge  to a constant in the $C^k$
sense, and consequently, by (4.7),  K\"ahler metrics
$\omega_{\tilde\phi_{t_i}}$ converge to a K\"ahler-Ricci soliton
$(\omega_g+\frac{\sqrt{-1}}{2\pi}\partial\overline\partial\overline\phi_{\infty})$
associated to the holomorphic vector field $X$. It remains to
prove that the limit
$(\omega_g+\frac{\sqrt{-1}}{2\pi}\partial\overline\partial\overline\phi_{\infty},
X)$  is  independent of the choice of sequence of
$\omega_{\tilde\phi_t}$. But the last  follows from the uniqueness
of K\"ahler-Ricci solitons proved in [TZ1] and [TZ2]. We leave the
details to reader. \qed\enddemo

\Refs\widestnumber\key{CTZ}

\item {[BS]} Batyrev, V. V. and Selivanova, E. N.,
      Einstein-K\"ahler metrics on symmetric toric Fano manifolds,
      J. Reine Angew. Math. 512, (1999), 225--236.

\item {[Ca]}Cao, H.D., Deformation of K\"ahler metrics to
K\"ahler-Einstein metrics on compact K\"ahler manifolds, Invent.
Math.,  81 (1985), 359-372.

\item {[CT1]}  Chen, X.X. and Tian, G., Ricci flow on
        K\"ahler-Einstein  surfaces, Invent. Math.,
        147 (2002), 487-544.

\item {[CT2]}  Chen, X.X. and Tian, G., Ricci flow on
        K\"ahler-Einstein  manifolds,  Duke Math. J., 131 (2006), 17-73.

\item {[CTZ]} Cao, H.D., Tian, G., and Zhu, X.H.,
      K\"ahler-Ricci solitons on compact K\"ahler manifolds with
      $c_1(M)>0$,  Geom. Anal. and Funct., 15 (2005), 697-719.

\item {[Ha]}  Hamilton, R.S., Three manifolds with positive Ricci
Curvature, J. Diff. Geom., 17 (1982), 255-306.

\item {[Ma]} Mabuchi, T., K-energy maps integrating Futaki
invariants, Toh\"oku Math. J., 38 (1986), 245-257.

\item {[Mi]}  Miguel, De G.,
       Differentiation of integrals in $\Bbb R^n$,
       Lectures in Math. 481 (1977), Springer-Verlag.

\item {[P1]} Perelman, G., The entropy formula for the Ricci flow
and its geometric applications, 2002, preprint.

 \item {[P2]} Perelman, G.,  unpublished.

\item {[ST]} Sesum, N. and Tian, G., Perelman's arguments for
uniform bounds on scalar curvature and diameter along the K\"ahler
Ricci flow, 2005, preprint.

\item {[TZ1]} Tian, G. and Zhu, X.H., Uniqueness of K\"ahler-Ricci
solitons, Acta Math., 184 (2000), 271-305.

\item {[TZ2]}  Tian, G. and Zhu, X.H.,
       A new holomorphic invariant and uniqueness of K\"ahler-Ricci
       solitons, Comm. Math. Helv.,  77 (2002),  297-325.

\item {[TZ3]}  Tian, G. and Zhu, X.H.,
       Convergence  of K\"ahler-Ricci flow, 2005,  to appear in Jour. of  Amer. Math. Soci..

\item{[WZ]} X-J, Wang and Zhu, X.H., K\"ahler-Ricci solitons on
toric manifolds with positive first Chern class., Advances in Math.,
188 (2004), 87-103.

\item {[Ya]} Yau, S.T.,
        On the Ricci curvature of a compact K\"ahler manifold and
        the complex Monge-Amp\`ere equation, I,
        Comm. Pure Appl. Math., 31 (1978), 339--411.

\item {[Zh]}   Zhu, X.H.,
       K\"ahler-Ricci soliton type equations on compact complex
       manifolds with $C_1(M)>0$,\newline
       J. Geom. Anal., 10 (2000),  759-774.

\endRefs
\enddocument